\documentclass[12pt,a4paper]{amsart}
\usepackage{url}
\usepackage{graphicx}
\usepackage{amssymb}
\usepackage{xy}
\xyoption{all}
\usepackage{psfrag}

\newcommand{\krho}{K^{\rho}}

\newcommand{\A}{\mathbb A}
\newcommand{\B}{\mathbb B}

\newcommand{\Z}{\mathbb Z}
\newcommand{\Q}{\mathbb Q}

\newcommand{\X}{\mathbb{X}}

\newcommand{\OO}{\mathcal O}

\newcommand{\dnorm}{\tilde{d}} 
\newcommand{\nuopen}{\buildrel _{\circ} \over {\nu}}  

\def\spinc{\ifmmode{\textrm{Spin}^c}\else{$\textrm{Spin}^c$}\fi}
\newcommand{\spincs}{\mathfrak s}

\def\oz{Ozsv{\'a}th-Szab{\'o}}

\newtheorem{theorem}{Theorem}[section]
\newtheorem{thm}{Theorem}

\newtheorem{cor}[thm]{Corollary}

\newtheorem{lemma}[theorem]{Lemma}

\newtheorem{proposition}[theorem]{Proposition}
\newtheorem{corollary}[theorem]{Corollary}
\theoremstyle{definition}
\newtheorem{definition}[theorem]{Definition}

\title{Concordance properties of parallel links}
\author[Daniel Ruberman]{Daniel Ruberman}
\address{Department of Mathematics, MS 050\newline\indent Brandeis
University \newline\indent Waltham, MA 02454}
\email{\rm{ruberman@brandeis.edu}}
\author[Sa\v so Strle]{Sa\v so Strle}
\address{Faculty of Mathematics and Physics \newline\indent 
University of Ljubljana \newline\indent Jadranska 21 \newline\indent 
1000 Ljubljana, Slovenia }
\email{saso.strle@fmf.uni-lj.si}
\renewcommand{\phi}{\varphi}
\thanks{The first author was partially supported by NSF Grant 0804760.   The second author was  supported in part by the Slovenian Research Agency program No. P1-0292-0101.  Visits of the authors were supported by a Slovenian-U.S.A. Research Project BI-US/09-12-004, and by NSF Grant 0813619.}
\begin{document}
\begin{abstract} We investigate the concordance properties of `parallel links' $P(K)$, given by the $(2,0)$ cable of a knot $K$. We focus on the question: if $P(K)$ is concordant to a split link, is $K$ necessarily slice?  We show that if $P(K)$ is smoothly concordant to a split link, then many smooth concordance invariants of $K$ must vanish, including the $\tau$ and $s$-invariants, and suitably normalized $d$-invariants of Dehn surgeries on $K$.   We also investigate the $(2,2\ell)$ cables $P_\ell(K)$, and find obstructions to smooth concordance to the sum of the $(2,2\ell)$ torus link and a split link.
\end{abstract}
\maketitle

\section{Introduction}
A central theme in the study of link theory has long been the relationship between the properties of a link and those of its individual components.   Recall that a split link is one whose components lie in disjoint balls; such a link is certainly determined by its individual components. It is 
classical~\cite{brunn,debrunner:brunnian} that not all links are split.  However, showing that 
not every link (say with linking number $0$) is concordant to a split link takes more work; in the classical dimension this was first done using Milnor's $\bar\mu$-invariants~\cite{milnor:mubar,milnor:link-groups}.  Examples in higher odd dimensions were first constructed by S.~Cappell and J.~Shaneson~\cite{cappell-shaneson:links} and independently by A.~Kawauchi~\cite{kawauchi:split}.    The even dimensional case remains stubbornly out of reach.

We investigate Kawauchi's construction in the classical dimension, where the additional distinction between the smooth and topological categories comes into play.  Kawauchi considers the $(2,0)$ cable of a knot $K$; we denote this `parallel' link by $P(K)$.  Kawauchi showed that if $P(K)$ is concordant to a split link, then $K$ is algebraically slice.  We show that if $P(K)$ is concordant to a split link, then a host of knot concordance invariants of $K$ must vanish.   Our results suggest the conjecture that $P(K)$ is smoothly (resp. topologically) concordant to a split link if and only if $K$ is smoothly (resp. topologically) slice.  In the last section, we consider the $(2,2\ell)$ cables, denoted $P_\ell(K)$, and give examples of knots $K$ for which $P_\ell(K)$ is topologically, but not smoothly concordant to the corresponding cable of the unknot.

Before stating our main theorem, we describe a convenient normalization of the Ozsv{\'a}th-Szab{\'o} $d$-invariant~\cite{oz:boundary} of surgery on an oriented knot.  Write $\spinc(Y)$ for the set of $\spinc$ structures on the manifold $Y$.  As we will detail below in Section~\ref{S:hf-prelim}, for any oriented knot 
$K$ in $S^3$, one can canonically label \spinc\ structures on $S^3_{p/q}(K)$ by elements $i\in \Z_p$.   We use the notation $\OO$ to denote the unknot.  
\begin{definition}\label{D:dnorm}
Let $K$ be a knot in $S^3$, and let  $i\in \Z_p$.  Then the normalized $d$-invariant $\dnorm(S^3_{p/q}(K),i) $ is defined to be 
$$
d(S^3_{p/q}(K),i) -d(S^3_{p/q}(\OO),i).
$$
\end{definition}
By Gordon's classic paper~\cite{gordon:contractible} and the homology cobordism invariance of the $d$-invariant~\cite{oz:boundary},  $d(S^3_{p/q}(K),i)$ (and hence $\dnorm(S^3_{p/q}(K),i)$) is a smooth knot concordance invariant.  (A recent preprint~\cite{peters:d-invariant} of T.~Peters considers the special case when the surgery coefficients are $\pm 1$.)  The normalization vanishes on the concordance class of the unknot and provides a neater statement of our main theorem. 
\begin{thm}\label{T:split-invariants} Suppose that $P(K)$ is smoothly concordant to a split link.  Then the following concordance invariants of $K$ must vanish:
\begin{enumerate}
\item $\tau(K)$, the Ozsv{\'a}th-Szab{\'o} $\tau$-invariant~\cite{oz:4ball}.  
\item $s(K)$, Rasmussen's invariant~\cite{rasmussen:khovanov-slice}.
\item $\delta_{2^n}(K)$, the $d$-invariant of the $2^n$-fold branched cover of $K$ \cite{manolescu-owens:delta}.
\item For any $p/q\in \Q-\{0\}$, and any $\spinc$ structure $i$ on $S^3_{p/q}(K)$, the normalized $d$-invariant $\dnorm(S^3_{p/q}(K),i)$.
\end{enumerate}
\end{thm}
The vanishing of each of the first three invariants follows directly from a simple geometric observation, Lemma~\ref{L:reflect}.  According to Tim Cochran, this same observation leads to the vanishing of invariants derived from the work of~\cite{cochran-orr-teichner:l2,cochran-orr-teichner:structure}, even if one only assumes a topological locally flat concordance to a split link. 
The last item in Theorem~\ref{T:split-invariants} is our main result, and makes use of some Dehn surgery manipulations, as well as Ozsv{\'a}th-Szab{\'o}'s recipe for computing~\cite{oz:q-surgery} the Heegaard-Floer homology of rational surgery on a knot. 
In fact, it follows from recent work of Ni and Wu~\cite{ni-wu:cosmetic} that the last item is equivalent to
\begin{itemize}
\item[($4'$)] $d(S^3_{\pm 1}(K))=0$;
\end{itemize}
see section 2.

Theorem~\ref{T:split-invariants}, combined with some known computations of $d$-invariants, leads directly to the following result. 
\begin{thm}\label{T:split} There are infinitely many concordance classes of two component links $L$ such that:
\begin{enumerate}
\item $L$ is topologically slice 
\item $L$ is a boundary link
\item $L$ is not smoothly  concordant to a split link.
\end{enumerate}
\end{thm}
Chuck Livingston has pointed out that links as described in Theorem~\ref{T:split} may also be obtained by Bing doubling a topologically slice knot.  We will explain the argument, which is basically that of Cimasoni~\cite{cimasoni:bing}, after the proof of Theorem~\ref{T:split}.  

It is a standard conjecture that the Bing double of $K$ is smoothly slice (which according to Lemma~\ref{L:split-link} below is equivalent to $B(K)$ being concordant to a split link) if and only if $K$ is.  The evidence for this is that many concordance invariants of $K$ vanish if $B(K)$ is slice~\cite{cha-livingston-ruberman:bing,cimasoni:bing,cochran-harvey-leidy:doubling, levine:doubles, vancott:bing}. We add one more to this list of invariants by showing the vanishing of the $\dnorm$-invariants of a knot $K$, given that $B(K)$ is smoothly slice.  This is deduced by a geometric argument (related to~\cite{cha-livingston-ruberman:bing}) establishing a link between concordance properties of $B(K)$ and those of $P(K)$. 
\begin{thm}\label{T:bing} 
Suppose that $B(K)$ is slice.  Then $P(K)$ is concordant, in a $\Z_2$-homology $S^3 \times I$, to a split link.
\end{thm}
Using Theorem~\ref{T:split}, this implies
\begin{cor}\label{C:bing-d}
If $B(K)$ is slice, then $\dnorm(S^3_{p/q}(K),i) = 0$ for any $p/q\in \Q-\{0\}$, and any $\spinc$ structure $i$ on $S^3_{p/q}(K)$. 
\end{cor}
With regard to the other parts of Theorem~\ref{T:split}, it was shown in~\cite{cha-livingston-ruberman:bing} that if $B(K)$ is slice, then $\tau(K) = \delta(K) = 0$ (the same argument would apply to $\delta_{2^n}(K)$).  On the other hand, it is not known if $s(K)$ would have to vanish.

In the last section we consider concordance of two component links with larger linking numbers.  A natural generalization of the question of whether $P(K)$ is concordant to a split link to the setting of links with nontrivial linking number $\ell$ is the question of whether $P_\ell(K)$ is concordant to $P_\ell(\OO)$ with knots tied in the individual components.   We will refer to that process as {\em local knotting} of a link.  Generalizing the above results we find links that are topologically concordant to the $(2,2\ell)$ torus link $P_\ell(\OO)$ but not smoothly concordant to a locally knotted $P_\ell(\OO)$. \\[2ex]
\textbf{Acknowledgments:}  We thank Jae Choon Cha, Tim Cochran, Matt Hedden, Chuck Livingston, Adam Levine, and Peter Ozsv\'ath for helpful discussions.

\section{Preliminaries}
All links will be assumed to be oriented. Generally speaking, concordance will refer to smooth concordance, with the adjective `topological' (always meaning locally flat) added as appropriate.    We will generally use the same letters for a link and its components, so that for example $L_1$ and $L_2$ would indicate the components of a two-component link $L$.

For a link $L=(L_1,\ldots,L_n)$ let $S(L)$ denote the split link with the same components as $L$, thus $S(L)=L_1\sqcup\cdots\sqcup L_n$, where the symbol `$\sqcup$' indicates that the components are in disjoint balls.  A local knotting of a link $J$ is any link obtained by tying knots in the components of $J$, or more formally by connected sum of $J$ with a split link.
We will make use of the following observation regarding concordance to split links; a more general version of this argument is given in~\cite[Lemma 2.5]{cha-kim-ruberman-strle:hopf}.

\begin{lemma}\label{L:split-link} 
If $L$ is concordant to a split link, then $L$ is concordant to $S(L)$.
\end{lemma}
\begin{proof}
Let $C$ be a concordance from $L$ to the split link $K_1 \sqcup \cdots \sqcup K_n$. Denote by $C_j$ the concordance from the component $L_j$ of $L$ to the corresponding component $K_j$.  Then $\tilde C$, defined by turning $C_1 \sqcup \cdots \sqcup C_n$ upside down, gives a concordance from $K_1 \sqcup \cdots \sqcup K_n$ to $S(L)$.  Composing $C$ with $\tilde C$ gives a concordance from $S(L)$ to $L$.
\end{proof}
\begin{figure}[h]
\psfrag{L1}{$L_1$}
\psfrag{L2}{$L_2$}
\psfrag{K1}{$K_1$}
\psfrag{K2}{$K_2$}
\psfrag{L}{$L$}
\psfrag{SL}{$S(L)$}
\psfrag{KK}{$K_1\sqcup K_2$}
\psfrag{C1}{$C_1$}
\psfrag{C2}{$C_2$}
\psfrag{C1t}{$\widetilde{C_1}$}
\psfrag{C2t}{$\widetilde{C_2}$}
\includegraphics[scale=.7]{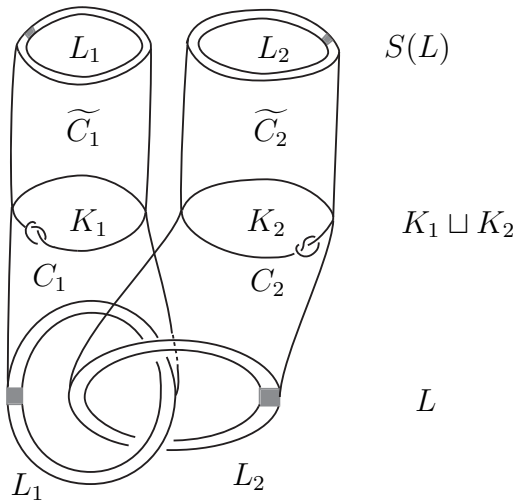}
\caption{Schematic illustration of proof of Lemma~\ref{L:split-link}}
\label{F:split-link}
\end{figure}

\subsection{Heegaard-Floer invariants}\label{S:hf-prelim}
We will make extensive use of the correction term, or $d$-invariant, introduced by \oz\ in~\cite{oz:4ball}. The $d$-invariant of a rational homology sphere $Y$ depends on the choice of $\spinc$ structure $\spincs$ on $Y$, and will be denoted $d(Y,\spincs)$.   Our main theorem is a comparison of the $d$-invariants of the Dehn surgery manifold $S^3_{p/q}(K)$ and those of the lens space $S^3_{p/q}(\OO) = -L(p,q)$.
We describe, briefly, a canonical way to enumerate \spinc\ structures on Dehn surgery on a knot in $S^3$, and hence a canonical correspondence between $\spinc(S^3_{p/q}(K))$ and $\spinc(S^3_{p/q}(\OO))$. The enumeration is in terms of relative \spinc\ structures on $S^3 - \nuopen(K) $, which by definition~\cite[Chapter I.4]{turaev:torsion-3mfd} are equivalence classes of non-singular vectors fields on $S^3 - \nuopen(K)$ that point outward at the boundary. These are determined by their relative first Chern classes in $H^2(S^3 - \nuopen(K),\partial\nu(K))$.  Hence~\cite{oz:knots} there is a one-to-one correspondence between the relative \spinc\ structures on $S^3 - \nuopen(K)$ and $S^3 - \nuopen(K')$ for any two oriented knots, where two relative $\spinc$ structures correspond if $c_1$ has the same evaluation on a Seifert surface.  This relative \spinc\ structure is labelled by $i \in \Z$ if the evaluation of the relative Chern class on a Seifert surface is $2i$.
Finally, a relative $\spinc$ structure 
determines a \spinc\ structure on $S^3_{p/q}(K)$, by extending a vector field on $S^3 - \nuopen(K)$ so that it is tangent to the core of the solid torus glued in by the Dehn surgery.  This \spinc\ structure is labelled by $i \pmod p$ in $\Z_p$.

The $d$-invariants of a general rational homology sphere $Y$ are difficult to compute, but there is a good deal known in the case when $Y$ is $p/q$ Dehn surgery on a knot $K$.  \oz~\cite{oz:q-surgery} give a  chain complex $\X_{i,p/q}$, described in terms of the knot chain complex $CFK^\infty(S^3,K)$, whose homology is $HF^+(S^3_{p/q}(K),i)$.  The chain complex $\X_{i,p/q}$ is the mapping cone of a map $D_{i,p/q}$ between chain complexes $\A_{i,p/q}$ and $\B_{i,p/q}$.  In turn, these are sums of subquotient complexes of $CFK^\infty(S^3,K)$, and the components of $D_{i,p/q}$ are defined in terms of certain maps $v_k$ and $h_k$ between those subcomplexes.

This description gives a good deal of information about $d(S^3_{p/q}(K),\spincs)$, and the implications have been elucidated in a recent preprint of Ni and Wu~\cite{ni-wu:cosmetic}.  The $U$-equivariance of $v_k$ and $h_k$ implies that they are determined, in sufficiently high gradings, by non-negative integers $V_k$ and $H_k$.  Moreover, $V_k$ are decreasing and vanish for $k \ge g$, $H_k$ are increasing and vanish for $k \le -g$ (where $g=g(K)$ denotes the genus of $K$), and $V_0=H_0$.  The $d$-invariants are then determined by these integers.
\begin{proposition}[{\cite[Proposition 2.11]{ni-wu:cosmetic}}]\label{P:ni-wu}
Suppose $p,\; q >0$, and $0\leq i \leq p-1$. Then
\begin{equation}\label{E:ni-wu}
\dnorm(S^3_{p/q}(K),i) = - 2 \max\{V_{\lfloor\frac{i}{q}\rfloor}, H_{\lfloor\frac{i-p}{q}\rfloor}\}.
\end{equation}
\end{proposition}

In view of the proposition the interesting range of numbers $V_k$ and $H_k$ is $V_i$ and $H_{-i}$ for $0\leq i \leq g-1$. This implies that the $d$-invariants of any large enough surgery on $K$ determine the $d$-invariants of all positive Dehn surgeries. The second item below was also observed by Ni and Wu.  
\begin{corollary}\label{C:dnorm}
{\rm(1)} Given knots $K,K'\subset S^3$, suppose that for some integer $n\ge 2\max\{g(K),g(K')\}-1$ the normalized $d$-invariants satisfy $$\dnorm(S^3_{n}(K),i)=\dnorm(S^3_{n}(K'),i)$$
for all $0\leq i \leq n-1.$ Then for any $p,q >0$ and all $0\leq i \leq p-1$ $$\dnorm(S^3_{p/q}(K),i)=\dnorm(S^3_{p/q}(K'),i).$$ 
{\rm(2)} If $\dnorm(S^3_1(K))=0$, then $\dnorm(S^3_{p/q}(K),i)=0$ for all $p,q >0$, and $0\leq i \leq p-1$.
\end{corollary}
\proof
If $n\ge 2g-1$, then for each $i$ satisfying $0\le i\le n-1$, at most one of $V_i$ or $H_{i-n}$ can be nonzero. Hence the normalized $d$-invariants determine all of these numbers, which in turn determine the normalized $d$-invariants of all positive Dehn surgeries.

If $\dnorm(S^3_1(K))=0$, then $V_0=H_0=0$ and all the relevant $V_k$ and $H_k$ vanish as well.
\endproof

\subsection{A Dehn surgery move}
In section~\ref{S:split}, we will use some moves on surgery diagrams with rational coefficients~\cite{rolfsen:knots}.  The first is the well-known {\em slam-dunk} move~\cite{cochran-gompf:donaldson,gompf-stipsicz:book}, pictured below.
\begin{figure}[h]
\psfrag{n}{$n$}
\psfrag{r}{$r$}
\psfrag{nr}{$n - \frac{1}{r}$}
\includegraphics[scale=.7]{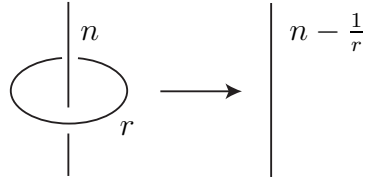}
\caption{Slam-dunk move; $n\in \Z$ and $r\in \Q$.}
\label{F:slam-dunk}
\end{figure}

\begin{figure}[ht]
\psfrag{K}{$K$}
\psfrag{K1}{$K_1$}
\psfrag{K2}{$K_2$}
\psfrag{pq}{$p/q$}
\psfrag{l}{$\ell$}
\psfrag{1q}{$\ell+1/q'$}
\psfrag{f}{${\displaystyle\ell+\frac{p}{q+ q'p}}$}
\includegraphics[scale=1.0]{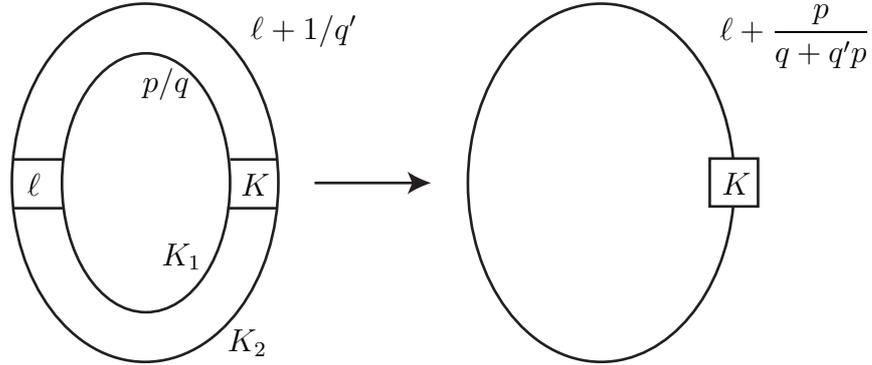}
\caption{Slap-shot move on $P_\ell(K)$}
\label{F:slap-shot}
\end{figure}

The second move is similar to the first one, and looks sort like a sideways slam-dunk.  In keeping with the sporting terminology, we call this the {\em slap-shot} move.
The move simplifies certain rational surgeries on the $(2,2\ell)$ cable of a knot $K$, which we label $P_\ell(K)$ (for related discussion see \cite{clark:cable}).  We do $p/q$ surgery on one copy, say $K_1$, and (for an integer $q'$) we do $\ell + 1/q'$ surgery on the other. The slap shot move asserts the diffeomorphism indicated in Figure \ref{F:slap-shot}.

The proof is quite similar to the proof of the slam-dunk.  First do the $\ell+1/q'$ surgery on $K_2$.  By an isotopy, $K_1$ may be moved inside the solid torus that is glued in during that first surgery; in fact, it will be isotopic to the core of that solid torus.  But then the $p/q$ surgery on $K_1$ turns that solid torus into another solid torus, which means that the surgery on $P_\ell(K)$ is diffeomorphic to a surgery on $K_1$.  It is easy to compute the surgery coefficient to be $\ell+p/(q+q'p)$.   One may alternatively establish the slap-shot move by a sequence of slam-dunks (and the inverse of the slam-dunk) and ordinary Rolfsen twists.

\section{Concordance to split links}\label{S:split}
In this section, we investigate the smooth concordance properties of the link $P(K)$ given by the $(2,0)$ cable of a knot $K$ in $S^3$.  Note that if $K$ is (topologically) slice, then $P(K)$ is (topologically) concordant to $P(\OO)$, and hence is slice.  Also, $P(K)$ bounds two parallel copies of a Seifert surface for $K$, and hence is a boundary link. As mentioned in the introduction, the work of Kawauchi~\cite[Theorem 5.1]{kawauchi:split} suggests the conjecture that if $P(K)$ is smoothly (resp.~topologically) concordant to a split link, then $K$ is smoothly (resp.~topologically) slice.  We show that if $P(K)$ is smoothly concordant to a split link, then various gauge-theoretic concordance invariants of $K$ must vanish.  Such results may be viewed as evidence for that conjecture.  

For a $k$-component link $L$, denote by $S^3_{r_1,\ldots,r_k}(L)$ the result of surgery on a link $L$, with surgery coefficients $r_1,\ldots,r_k \in \Q$.  Likewise, we indicate by $M_n(L)$ 
the $n$-fold cyclic branched cover of $S^3$, branched along all components of $L$, corresponding to the homomorphism $\pi_1(S^3 - L) \to \Z_n$ taking all meridians to $1\in \Z_n$.  For an oriented knot $K$, the notation $K^r$ indicates the same knot with reversed orientation, and $\krho$ the image of $K$ under a reflection of $S^3$.  

The first observation about $P(K)$ is close to saying that $K$ has order $2$ in the concordance group.
\begin{lemma}\label{L:reflect}
If $P(K)$ is concordant to a split link, then $K$ is concordant to $\krho$.   The oriented meridian of $\krho$ is homologous to the oriented meridian of $K$ in the complement of the concordance.
\end{lemma}
\begin{proof}
By Lemma~\ref{L:split-link}, if $P(K)$ were concordant to a split link, it would be concordant to $K\sqcup K$.  Orient $P(K)$ so that the components have opposite orientation, which implies that those components cobound an oriented annulus.  Glue that annulus to a cobordism between $P(K)$ and $K\sqcup K$ to obtain an annulus in $B^4$ with oriented boundary $K\sqcup K^r$.  The knots $K$ and $K^r$ lie in disjoint $3$-balls in $S^3$; viewing $B^4$ as a product of one of those with an interval, the annulus becomes a cobordism between $K$ and $(K^r)^{\rho r} = \krho$.
\end{proof}
The reader who is confused about the two sorts of orientations should compare this argument with the standard proof that the inverse of $K$ in the concordance group is given by $-K = K^{\rho r}$.  

Lemma~\ref{L:reflect} gives rise to additional restrictions to $P(K)$ being concordant to a split link, via any $\Z$--valued homomorphism from the concordance group that is insensitive to the orientation on $K$.  So for instance, we immediately deduce the first 3 parts of Theorem~\ref{T:split-invariants}.
\begin{corollary}\label{C:split-tau-s}
If $P(K)$ is concordant to a split link, then $\tau(K) = s(K) =\delta_{2^k}(K)=0$, where $s$ is Rasmussen's invariant~\cite{rasmussen:khovanov-slice}, and $\delta_{2^k}(K)$ is the $d$-invariant for the spin structure on the $2^k$-fold cover of $S^3$ branched along $K$ (cf.~\cite{manolescu-owens:delta}).
\end{corollary}

\begin{proof}[Proof of Theorem~\ref{T:split}]
Consider the parallel $P(K)$, where $K = D^+(J,0)$ is the untwisted, positive-clasped Whitehead 
double of $J$.  The algorithm of Akbulut and Kirby~\cite{akbulut-kirby:branch} exhibits the double branched cover $M_2(K)$ as surgery on a link of two 
components, one of which is $J \# J^r$ with framing $0$, and the other of which is a meridian of that knot, 
with framing $-2$.  An application of the slam-dunk move then identifies $M_2(K)$ with $S^3_{1/2}(J\#J^r)$, hence if $P(K)$ is concordant to a split link, then by Corollary \ref{C:split-tau-s}
$d(S^3_{1/2}(J\#J^r))=0$.

At this point, there are many choices for the knot $J$ that will lead to a contradiction.  For example, we may take $J$ to be the connected sum of $n$ copies of the right-handed trefoil.  Alternatively, we choose $J$ to be the $(2,2n+1)$ torus knot. From this point the argument proceeds as in the proof of Theorem 4.1 of~\cite{cha-kim-ruberman-strle:hopf}. In particular, the $d$-invariant of $S^3_{1/2}(J\#J^r)$ is $-2n$ from which it also follows that different choices for $K$ are not concordant and hence neither are the corresponding links.
\end{proof}

Chuck Livingston pointed out an alternate construction that proves Theorem~\ref{T:split}, with the additional feature that the components are unknotted.   As observed in~\cite[Propositions 1.1 and 3.2]{cimasoni:bing} Bing doubles are always boundary links, and there is a genus-$0$ cobordism between the Bing double of a knot $K$ and the untwisted positive Whitehead double of $K$.  This proves that the Bing doubles of many knots (i.e. those whose Whitehead doubles are not smoothly slice) are not smoothly slice.  By Lemma~\ref{L:split-link}, such a Bing double will not be concordant to a split link, because its components are unknotted.

Another collection of gauge-theoretic concordance invariants of a knot $K$ is given by the set of 
$d$-invariants~\cite{oz:boundary} of all of the Dehn surgeries on $K$.  These have been investigated for surgery coefficients $\pm 1$ by T.~Peters~\cite{peters:d-invariant}, but it is not hard to show that the concordance invariance extends to any surgery coefficient and any $\spinc$ structure.   In seeking to show that a knot $K$ is not smoothly slice, one would want to show that the $d$-invariants of $S^3_{p/q}(K)$ differ from those of $p/q$ surgery on the unknot (for corresponding $\spinc$ structures).  

Hence we would like to show that if $P(K)$ is concordant to a split link, then the $d$-invariants of $S^3_{p/q}(K)$ are equal to those of  $S^3_{p/q}(\OO)$. 
We start by deducing some homology cobordisms from a concordance of $P(K)$ to a split link.
\begin{proposition}\label{P:split}
If $P(K)$ is concordant to a split link, then for any $p/q\in \Q$ and $q'\in\Z$
\begin{enumerate}
\item $S^3_{-p/q}(K) $ is homology cobordant to $-S^3_{p/q}(K)$, and
\item $S^3_{p/q}(K) \# S^3_{1/q'}(K)$  is homology cobordant to $S^3_{p/(q+q'p)}(K)$.
\end{enumerate}
\end{proposition}
\begin{proof}
In general, if a link $L$ is concordant to a split link there is a homology cobordism between $S^3_{r_1,\ldots,r_k}(L)$ and $S^3_{r_1}(L_1) \# \cdots \# S^3_{r_k}(L_k)$.  This follows by Lemma~\ref{L:split-link}, together with Gordon's paper~\cite{gordon:contractible}.  
We use the orientation-preserving diffeomorphism $S^3_{-p/q}(K) \cong  - S^3_{p/q}(\krho)$ for any $p/q \in \Q \cup \{\infty\}$ combined with Lemma~\ref{L:reflect} to get the homology cobordism of item (1).

The proof of the second item again starts by surgering a concordance between $P(K)$  and 
$K\sqcup K$ to get a homology cobordism between surgery on $S^3_{p/q,1/q'}(P(K))$ and $S^3_{p/q}(K) \# S^3_{1/q'}(K)$.  Then an application of the slap-shot move to $S^3_{p/q,1/q'}(P(K))$ gives (2).
\end{proof}

Proposition~\ref{P:split} leads to obstructions, by using the  Ozsv{\'a}th-Szab{\'o} $d$-invariant \cite{oz:boundary}.  
We note that a homology cobordism induces a canonical bijection on the set of $\spinc$ structures. In fact, in the case of interest, the orientation of a knot induces an orientation of any knot to which it is concordant, hence the labelling of relative $\spinc$ structures of the two knots agrees.
\begin{proposition}\label{P:split-d}
If $L = P(K)$ is concordant to a split link, then for any $p/q\in \Q-\{0\}$ and any $q'\in\Z$, 
\begin{enumerate}
\item $d(S^3_{-p/q}(K),i) =-d(S^3_{p/q}(K), i)$, and
\item $d(S^3_{p/q}(K),i)+ d( S^3_{1/q'}(K)) = d(S^3_{p/(q+q'p)}(K),i)$
\end{enumerate}
for all $\spinc$ structures  $i \in\Z_p$.
\end{proposition}
\begin{proof} These follow from the homology cobordism invariance of the $d$-invariant.\end{proof}
We note that the $d$-invariants of the $0$ surgery on $K$ are determined by the above (see \cite[Proposition 4.12]{oz:boundary}).

Proposition~\ref{P:split-d} can be used by itself to obstruct $P(K)$ being concordant to a split link. It gains strength by being coupled with relations between $d$-invariants for different surgeries on an arbitrary knot.   This is how we prove part $(4)$ of Theorem~\ref{T:split-invariants}.

\begin{theorem}\label{T:vanishing}
If $P(K)$ is concordant to a split link, then 
$$d(S^3_{\pm 1}(K))= 0. $$ 
Moreover, for any $p/q\in \Q-\{0\}$, and any $\spinc$ structure $i$ on $S^3_{p/q}(K)$, the normalized $d$-invariant $\dnorm(S^3_{p/q}(K),i)$ vanishes.
\end{theorem}
\begin{proof}
Part (2) of Proposition~\ref{P:split-d} for $p=q=q'=1$ yields 
$$2d(S^3_{1}(K))=d(S^3_{1/2}(K)).$$ 
Combining this with $d(S^3_{1}(K))=d(S^3_{1/2}(K))$ (see \cite[Corollary 9.14]{oz:boundary}) we obtain the result for $1$ surgery on $K$. Vanishing of other normalized $d$-invariants now follows from part (2) of Corollary \ref{C:dnorm} and part (1) of Proposition~\ref{P:split-d}.
\end{proof}

\section{Bing doubles}
The relationship between $B(K)$ and Kawauchi's results on cables of $K$ is a key ingredient in~\cite{cha-livingston-ruberman:bing}.  In this section we develop this parallel a bit more, yielding new restrictions on a  knot whose Bing double is slice.
\begin{proof}[Proof of Theorem~\ref{T:bing}]
We make use of the covering link technique in~\cite[\S3]{cha-livingston-ruberman:bing}. Supposing that $B(K)$ is slice, consider the lift of one component of $B(K)$ to the $4$-fold branched cyclic cover over the other component, drawn in~\cite[Figure 3]{cha-livingston-ruberman:bing}.  By lifting the corresponding component of the concordance, we obtain a concordance in a $\Z_2$-homology $S^3 \times I$ (say $W$) between the $2$-component link $(J_1, J_2)$ in that figure (consisting of two `adjacent' lifts)
and the unlink. Direct inspection of that figure shows that $(J_1, J_2)$ is, in our notation, $P(K) \# (K^r \sqcup K^r)$.

Let $C$ be a concordance in $W$ between $P(K) \# (K^r \sqcup K^r)$ and the unlink, and consider embedded arcs $\alpha_1$, $\alpha_2$ in the two components of $C$.  Take the connected sum of $C$, along $\alpha_1$, $\alpha_2$, with two copies of the product concordance between $K$ and itself.  The result is then a concordance, in $W$, between $P(K) \# ((K^r \# K) \sqcup (K^r \# K))$ and the split link $K \sqcup K$.
On the other hand, $P(K) \# ((K^r \# K) \sqcup (K^r \# K))$ is concordant to $P(K)$, since $K^r \# K$ is slice.  Composing the two concordances proves the theorem.
\end{proof}
For the proof of Corollary~\ref{C:bing-d}, we need only observe that the proofs of the vanishing of the $\dnorm$-invariants given in Section~\ref{S:split} require only that $P(K)$ be concordant to a split link in a $\Z_2$-homology $S^3 \times I$, both of whose boundary components are $S^3$.

\section{Other linking numbers}
For any integer $\ell$ and knot $K$ denote by $P_\ell(K)$ the $(2,2\ell)$ cable of $K$; note that the parallel $P(K)$ is equal to $P_0(K)$. In this section we consider the question whether $P_\ell(K)$ is concordant to a locally knotted $P_\ell(\OO)$. In a similar vein as for $\ell=0$ we give examples of links for which this holds in the topological category but not in the smooth.  We start with two simple results analogous to the first two parts of Theorem~\ref{T:split-invariants} that yield single examples where there is a topological concordance to a locally knotted $P_\ell(\OO)$, but no such smooth concordance.  The main result of this section is Theorem~\ref{T:linking}, which gives infinitely many examples. We conjecture that if $P_\ell(K)$ is concordant to a locally knotted $P_\ell(\OO)$ then $K$ must be a slice knot. 

\begin{proposition}\label{P:twisted-tau}
Suppose that $P_\ell(K)$ is concordant to a locally knotted $P_\ell(\OO)$.  Then, for any $\ell$, both  $\tau(K)$ and $s(K)$ must vanish.
\end{proposition}
\begin{proof}
According to~\cite[Lemma 2.5]{cha-kim-ruberman-strle:hopf} if $P_\ell(K)$ is concordant to a locally knotted $P_\ell(\OO)$, then we may assume that this locally knotted link is $P_\ell(\OO)\#(K \sqcup K)$.  Orient $P_\ell(K)$ so that the components go in opposite directions, or in other words so that $P_\ell(K)$ bounds an annulus.  It follows that $P_\ell(\OO)\#(K \sqcup K)$ bounds an annulus in $B^4$.  

Do a band sum of the two components of $P_\ell(\OO)\#(K \sqcup K)$ indicated in Figure~\ref{F:l-split} to obtain $K \# K^r$.  \begin{figure}[ht]
\psfrag{K}{$K$}
\psfrag{l}{$\ell$}
\includegraphics[scale=1.0]{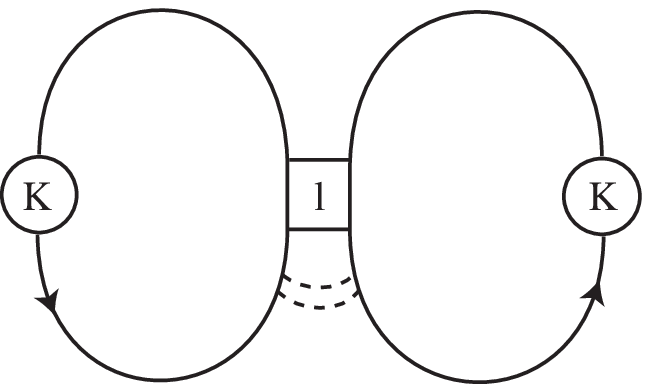}
\label{F:l-split}
\end{figure}
We conclude that $K \# K^r$ is the boundary of a genus-$1$ surface in $B^4$. From the adjunction inequalities satisfied by the two invariants~\cite{oz:4ball,rasmussen:khovanov-slice}, we get that 
\begin{align*}
2|\tau(K)| &= |\tau(K \# K^r)| \leq 1,\quad\text{and}\\
2|s(K)| &= |s(K \# K^r)| \leq 1.
\end{align*}
Since $s(K)$ and $\tau(K)$ are both integers, they must vanish.
\end{proof}

Choosing $K$ to be the untwisted positive Whitehead double of a knot with $\tau <0$, we conclude from~\cite{hedden:double} that while $P_\ell(K)$ is topologically concordant to $P_\ell(\OO)$, it is not smoothly concordant to a locally knotted $P_\ell(\OO)$.   To get infinitely many concordance classes of such examples we make use of the $d$-invariants, which have a greater range than the $\Z$-valued $s$ and $\tau$-invariants.  We begin by establishing a useful way to estimate the $d$-invariants of Dehn surgery on a knot.

\begin{lemma}\label{L:d-bounded}
For any nonnegative integers $n,p$ and any $r\in\Q$ there is a constant $d(n,p,r)$ such that for any knot $K$ that can be unknotted by changing $p$ positive and $n$ negative crossings the correction terms of the $r$ surgery on $K$ satisfy
\begin{align*}
&d(S^3_r(K),\spincs)\ge -d(n,p,r) \ \text{for}\ r\not\in[0,4p]\\
&d(S^3_r(K),\spincs)\le d(n,p,r) \ \text{for}\ r\not\in[-4n,0]
\end{align*}
for all $\spinc$ structures $\spincs$.
\end{lemma}
\begin{proof}
Suppose $K$ can be unknotted by changing $p$ positive and $n$ negative crossings. Realizing the crossing changes by surgeries on $-1$ framed unknots linking $K$ at the appropriate crossings yields the unknot in $S^3$ with framing $r+4n$. The cobordism $W$ corresponding to the surgeries is negative definite for $r\not\in[-4n,0]$. To see this note that it is sufficient to treat one crossing change at a time.

The cobordism resulting from changing a single crossing is given by adding a $2$-handle to $S_r^3(K)$ along the $J$ pictured in blue in Figure~\ref{F:crossing}; its framing in $S^3$ is $-1$.  If the crossing is positive, then $J$ is null-homologous in $S_r^3(K)$, and the generator of $H_2(W)$ is represented by a null-homology for $J$ in $S^3 -K$  together with the core of the two-handle, and has self-intersection $-1$.  So the cobordism is always negative definite in this case.

In case of a negative crossing consider the cobordism $W'$ corresponding to the surgeries on both $K$ and $J$. Suppose first $r>0$ and expand $r=a/b$ into a continued fraction. Using the slam-dunk move this gives a presentation of the $r$ surgery on $K$ as an integral surgery on a link. Denote the intersection form of the resulting positive definite cobordism by $A$ and the minor of $A$ obtained by deleting the generator corresponding to $K$ by $B$. Then $\det A=a$ and $\det B=b$. The determinant of the intersection pairing on $W'$ is then $-a-4b=-b(r+4)$ which is negative. This proves that the cobordism $W$ is negative definite. If $r<0$ write $-r=a/b$. In this case $A$ is negative definite, $\det A=(-1)^k a$, $\det B=(-1)^{k-1}b$ and the determinant of the intersection pairing on $W'$ is $(-1)^k b(r+4)$. This has the opposite sign as $\det A$ for $r<-4$ showing again that for such $r$ the cobordism $W$ is negative definite.

Alternatively, we can identify the generator of $H_2(W)$ in the case of negative crossing as well. The situation is a bit more complicated, as the self-intersection of the generator is a rational number.  Assume $p$ is odd; if $p$ is even then replace $p$ by $p/2$ in what follows. The generator is represented by a null-homology for $p$ times $J$ in $S_r^3(K)$, union $p$ copies of the core of the two handle. The self-intersection of the generator is a rational number, given by $1/p$ times the intersection of a push-off of $J$ (using the given framing) with such a null-homology. We may compute this via a formula of Hoste~\cite{hoste:casson} (he treats the setting where $r = \pm 1$, but the method works in general) and get $-1 - 4/r$.  The sign of this is negative for $r \not\in [-4,0]$, so the cobordism is negative definite except for those values of $r$.

\begin{figure}
\psfrag{J}{$J$}
\psfrag{K}{$K$}
\psfrag{r}{$r$}
\psfrag{rp4}{$r+4$}
\psfrag{m1}{$-1$}
\psfrag{neg}{negative}
\psfrag{pos}{positive}
\includegraphics[scale=.9]{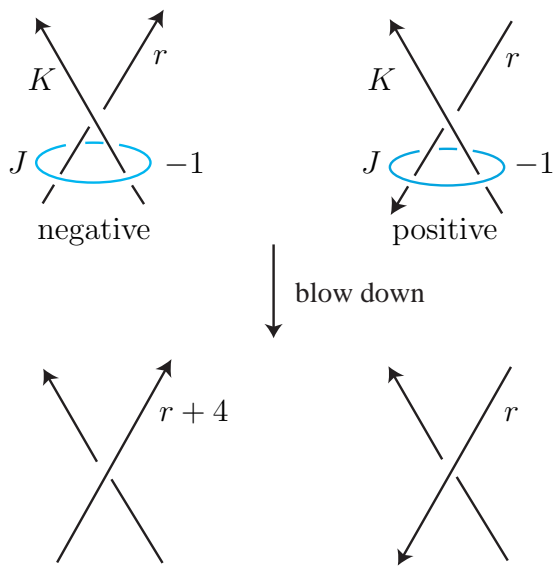}
\caption{Changing crossings}
\label{F:crossing}
\end{figure}
Then by \cite[Theorem 9.6]{oz:boundary} for any $\spinc$ structure $\spincs$ on the cobordism we have
$$4d(S^3_r(K),\spincs')+c_1(\spincs)^2+p+n\le 4d(S^3_{r+4n}(\OO),\spincs''),$$
where the primes denote the restrictions of $\spinc$ structures to the boundary components. Note that the restriction map $H^2(W;\Z)\to H^2(S^3_r(K);\Z)$ is onto since $H^3(W,S^3_r(K);\Z)\cong H_1(W,S^3_{r+4n}(\OO);\Z)=0$ hence the above inequality yields the desired upper bounds independent of $K$.

To obtain the lower bounds one could consider the positive definite cobordism corresponding to unknotting $K$ by $+1$ surgeries and then changing the orientation. Alternatively we note that reversing the roles of the unknot $\OO$ and $K$ one can obtain $K$ from an appropriate diagram for $\OO$ by changing $n$ positive and $p$ negative crossings. This gives a negative definite cobordism from $S^3_{r-4p}(\OO)$ to $S^3_r(K)$ for $r\not\in [0,4p]$.
\end{proof}

\begin{theorem}\label{T:linking}
For any integer $\ell$ there are infinitely many concordance classes of two component links $L$ that are topologically concordant to $P_\ell(\OO)$ but not smoothly concordant to a locally knotted $P_\ell(\OO)$.
\end{theorem}
\begin{proof}
Consider the links $P_\ell(K(n))$ for $K(n)=D^+(J(n),0)$ the positive Whitehead double of the $(2,2n+1)$ torus knot $J(n)$. Orient $P_\ell(K(n))$ so that the components have opposite orientation and choose $\ell<0$ so that the linking number of $P_\ell(K(n))$ is $|\ell|$.  (The case when $\ell>0$ is handled by reflecting the link.) The branched double cover $M_2(n)$ of $P_\ell(K(n))$ is $S^3_{2\ell}(K(n)\# K(n)^r)$. Since $K(n)$ can be unknotted by changing one positive crossing it follows from Lemma \ref{L:d-bounded} that the $d$-invariants of $M_2(n)$ are bounded below independent of $n$. Suppose $P_\ell(K(n))$ were concordant to a locally knotted $P_\ell(\OO)$.  As above, we may assume that that this locally knotted link is $P_\ell(\OO)\#(K(n) \sqcup K(n))$, whose branched double cover $M_2(n)'$ is $L(2\ell,1)\# 2M_2(K(n))$. The branched double cover $M_2(K(n))$ of $K(n)$ is $S^3_{1/2}(J(n) \# J(n)^r)$ whose $d$-invariant is $-2n$ (see proof of~\cite[Theorem 4.1]{cha-kim-ruberman-strle:hopf}) which yields a contradiction for large $n$.
\end{proof}

\def\cprime{$'$}
\providecommand{\bysame}{\leavevmode\hbox to3em{\hrulefill}\thinspace}

\end{document}